\documentclass[12pt]{article}
\usepackage{amssymb,amsfonts} %
\usepackage{amsthm} %
\usepackage{color}
\definecolor{darkgreen}{rgb}{0,.5,0}
\definecolor{darkred}{rgb}{.5,0,0}
\definecolor{darkcyan}{rgb}{0,.5,.5}
\definecolor{darkmagenta}{rgb}{.5,0,.5}
\usepackage[dvips,
            colorlinks,%
            linkcolor=darkred,%
            pagecolor=red,%
            filecolor=darkcyan,%
            citecolor=darkgreen,%
            urlcolor=darkmagenta]{hyperref} 
\def\proofr{\begin{proof}}
\def\proofend{\end{proof}}

\theoremstyle{plain} 
\newtheorem{theorem}{Theorem}[section]
\newtheorem{lemma}[theorem]{Lemma}
\newtheorem{proposition}[theorem]{Proposition}
\newtheorem{corollary}[theorem]{Corollary}
\theoremstyle{definition} 
\newtheorem{example}[theorem]{Example}
\newtheorem{remark}[theorem]{Remark}

\newcommand{\bpro}{\begin{proposition}}
\newcommand{\epro}{\end{proposition}}

\newcommand{\bcor}{\begin{corollary}}
\newcommand{\ecor}{\end{corollary}}

\newcommand{\bnote}{\begin{remark}}
\newcommand{\enote}{\end{remark}}

\newcommand{\btheorem}{\begin{theorem}}
\newcommand{\etheorem}{\end{theorem}}

\newcommand{\blemma}{\begin{lemma}}
\newcommand{\elemma}{\end{lemma}}

\newcommand{\bexam}{\begin{example}}
\newcommand{\eexam}{\end{example}}

\begin{document}
\title{On decomposability of $4$-ary distance $2$ MDS codes, double-codes, and
$n$-qua\-si\-groups of order $4$%
\thanks{Some results of this paper were presented at
the 8th International Workshop on Algebraic and Combinatorial Coding Theory
ACCT-VIII, which was held in Tsarskoe Selo (Russia, June 2002).}
}

\author{Denis~S.~Krotov%
\thanks{Sobolev Institute of Mathematics, Ak. Koptyuga, 4,
Novosibirsk, 630090, Russia (e-mail: \texttt{krotov@math.nsc.ru})}
}

\maketitle
\begin{abstract}
A subset $S$ of $\{0,1,\ldots,2t-1\}^n$
is called a $t$-fold MDS code if
every line in each of $n$ base directions contains exactly $t$ elements of $S$.
The adjacency graph of a $t$-fold MDS code is not connected
if and only if the characteristic function of the code is the repetition-free
sum of the characteristic functions of $t$-fold MDS codes of smaller lengths.

In the case $t=2$, the theory has the following application.
The union of two disjoint $(n,4^{n-1},2)$
MDS codes in $\{0,1,2,3\}^n$ is a double-MDS-code. If the adjacency
graph of the double-MDS-code is not connected, then the double-code
can be decomposed into double-MDS-codes of smaller lengths.
If the graph has more than two connected components, then the
MDS codes are also decomposable. The result has an interpretation
as a test for reducibility of $n$-qua\-si\-groups of order $4$.

Keywords:
\ifx\undefined\sep \def\sep{, } \fi
\ifx\undefined\MSC \def\MSC{MSC: } \fi
MDS codes\sep
$n$-quasigroups\sep
decomposability\sep
reducibility\sep
frequency hypercubes\sep
latin hypercubes

\MSC
05B99\sep 20N15\sep 94B25
\end{abstract}

\section{Introduction}

We consider the subsets $S$ of $\{0,1,\ldots,q-1\}^n$, where $q\geq 4$ is even,
with the following property:
every line in each of $n$ base directions contains exactly $q/2$ elements of $S$.
We call such objects $q/2$-fold MDS codes.
This paper establish a connection between the connectivity of a $q/2$-fold MDS code
and its decomposability.
More accurately, by the example of $q=4$,
we prove that the adjacency graph of a $q/2$-fold MDS code is not connected
if and only if the characteristic function of the code is the repetition-free
sum of the characteristic functions of $q/2$-fold MDS codes of smaller lengths.

$q/2$-Fold MDS codes are very natural objects of study;
they can be considered as a partial case of (strongly defined) frequency hypercubes,
an $n$-dimensional generalization of frequency squares
(questions of connectivity for a partial type of frequency squares were considered
in \cite{HiltonOth2003}).
Nevertheless, our research is motivated by studying
the $4$-ary distance $2$ MDS codes or, equivalently,
the $n$-quasigroups of order $4$.

A distance $2$ MDS code is decomposable if it can be represented
as a ``concatenation'' (see (\ref{e3}) in Section~\ref{sect:MDS}) of MDS codes of smaller length.
The goal of this work is to prove the following test for decomposability of
$4$-ary distance $2$ MDS codes.

Let $C$ and $C'$ be two disjoint MDS codes in $\{0,1,2,3\}^n$.
Assume that the adjacency graph of their union ($2$-fold MDS code)
has more than the minimal ($1$ or $2$)
and less than the maximal ($2^{n-1}$) number of connected components.
Then the MDS codes $C$ and $C'$ are decomposable.
Note that if $C'=\pi C$ where $\pi$  is a permutation of type $(a,b)(c,d)$
of the alphabet symbols in one coordinate,
then there are at least $2$ connected components.
Otherwise the minimal number of components is $1$.
If the adjacency graph of $C\cup C'$ has $2^{n-1}$ connected components, then
$C$ and $C'$ belong to the class of semilinear MDS codes
(see Section~\ref{sect:MDS}).

In particular, this test means that we cannot get a ``new'' code
if we combine parts of two disjoint
$4$-ary distance $2$ MDS codes $C_1$ and $C_2$
(see Theorem~\ref{th:switch} for the details).
So, this ``switching'' method, which works well, for example, for constructing
$1$-perfect binary codes with nontrivial properties (see e.\,g. \cite{Sol:switchings}),
cannot provide something interesting in the case of $4$-ary distance $2$ MDS codes.

Since there is a one-to-one correspondence between $q$-ary distance $2$ MDS codes
of length $n+1$ and $n$-qua\-si\-groups of order $q$
(the value arrays of $n$-qua\-si\-groups are also known as latin $n$-cubes,
an $n$-dimensional
generalization of latin squares),
we can interpret the results in terms of
$n$-qua\-si\-groups of order $4$ (Section~\ref{sect:quasigroup}).
The decomposability, or reducibility, of $n$-qua\-si\-groups is a natural concept;
for arbitrary order it was considered, for example, in \cite{Belousov,Cher}.

The mention of $1$-perfect binary codes above is not an accident.
There are concatenation constructions of such codes \cite{Phelps84,ZinLob:2000}
based on distance $2$ MDS codes, or $n$-qua\-si\-groups.
Moreover, as shown in \cite{AvgHedSol:class},
any $1$-perfect binary code of length $m$ and rank $\leq \mbox{minimal\_rank}{+} 2$ 
is described by a collection of 
distance $2$ MDS codes in $\{0,1,2,3\}^{(m+1)/4}$
(the rank is the dimension of the code linear span; the minimal rank
of a $1$-perfect binary code is $m-\log_2 (m+1)$).
So, the properties of distance $2$ MDS codes are closely related to properties
of some $1$-perfect codes.

Concepts closely related with $4$-ary distance $2$ MDS codes are the
concepts of a double-code and a double-MDS-code
(i.\,e., $2$-fold MDS codes in $\{0,1,2,3\}^n$).
Double-codes and double-MDS-codes
have many useful properties, which are discussed in Section~\ref{sect:pre}.
Studying MDS codes, we can think that a double-MDS-code is
the union of two disjoint
$4$-ary distance $2$ MDS codes and a double-code is
a part of a double-MDS-code closed with respect to adjacency.
In fact, there are double-MDS-codes, as well as $q/2$-fold MDS codes,
that are not splittable into distance $2$ MDS codes,
see \cite{KroPot:nonsplittable},
and the class of all double-MDS-codes can be considered independently.

In Section~\ref{sect:def} we give main definitions and notations.
In particular, we define the concept of a double-MDS-code,
which is a set with properties of the union of two disjoint distance $2$ MDS codes.
In Section~\ref{sect:pre} we prove some preliminary results.
In Section~\ref{sect:2MDS} we prove the theorem on
the decomposition of double-MDS-codes into prime double-MDS-codes
and show how to generalize the result to $q/2$-fold MDS codes.
In Sections~\ref{sect:MDS} and~\ref{sect:quasigroup} we discuss
the decomposability of distance $2$ MDS codes and $n$-qua\-si\-groups.
In the Appendix we prove some auxiliary lemmas about functions with separable arguments,
which are used in Sections~\ref{sect:pre} and~\ref{sect:2MDS}.

The author would like to thank Vladimir Potapov and Sergey Avgustinovich
for many helpful discussions and suggestions.

\section{Basic notations and definitions} \label{sect:def}

Let ${\Sigma}\triangleq\{0,1,2,3\}$ and ${\Sigma}^n$ be the set of words of
length $n$ over the alphabet ${\Sigma}$.
Denote $[n] \triangleq \{1,\ldots,n\}$.
For $\bar x=(x_1,x_2,\ldots ,x_n)$ we
use the following notation:
\def\laaa{{[}}
\def\raaa{{]}}
\def\lbbb{{[}}
\def\rbbb{{]}}
\def\lbb{{[}}
\def\rbb{{]}}
\begin{eqnarray*}
\bar x^{\laaa k \raaa}\lbb y\rbb &\triangleq& (x_1,\ldots ,x_{k-1},y,x_{k+1},\ldots ,x_n),\\
\bar x^{\laaa k_1,k_2,\ldots ,k_s\raaa}\lbbb y_1,y_2,\ldots ,y_s \rbbb &\triangleq&
\bar x^{\laaa k_1\raaa}\lbbb y_1\rbbb ^{\laaa k_2\raaa}\lbbb y_2\rbbb \ldots \vphantom{|}^{\laaa k_s\raaa}\lbbb y_s\rbbb.
\end{eqnarray*}

A set of four elements of ${\Sigma}^n$ that differ in only one (\,$i$th\,) coordinate
is called a \em line \em ($i$-\em line\em) of ${\Sigma}^n$.
Let ${\mathcal E}_i(\bar x)$ denote the $i$-line that contains $\bar x\in {\Sigma}^n$.
If $S\subset {\Sigma}^n$, then
\[
{\mathcal E}_i(S)\triangleq\bigcup_{\bar x\in S}{\mathcal E}_i(\bar x)
\]
(the union of the $i$-lines through the points of $S$) and
\[
{\mathcal F}_{i,j;\bar x}
S\triangleq\{(b,c)\in {\Sigma}^2:\bar x^{\laaa i,j\raaa}\lbbb b,c\rbbb\in S\}
\]
(the cut of $S$ in the ``$i,j$-plane'' through $\bar x$).

A set $C\subset {\Sigma}^n$ is called a \em $4$-ary distance $2$  MDS
code $($of length {$n$}$)$ \em or \em $(n,2)_4$ MDS code \em if each
line of ${\Sigma}^n$ contains exactly one element of $C$. A function
$g:{\Sigma}^n\to {\Sigma}$ is called an {\it $n$-qua\-si\-group of order $4$} if for
each $i\in [n]$ and $y,x_1,\ldots ,x_{i-1},x_{i+1},\ldots ,x_n\in {\Sigma}$ there exists
$x_i=g^{\langle i\rangle}(x_1,\ldots ,x_{i-1},y,x_{i+1},\ldots ,\linebreak[1]x_n)\in {\Sigma}$
such that $y=g(x_1,\ldots ,x_n)$.
Clearly, the function $g^{\langle i\rangle}$ is also an $n$-qua\-si\-group of order $4$.
For the rest of the paper we omit the words
``$4$-ary distance $2$'' and ``of order $4$'' because we consider
only MDS codes and $n$-qua\-si\-groups with such parameters.
The following one-to-one correspondence between MDS codes and
$n$-qua\-si\-groups is obvious and well known.
\bpro\label{MDS<->quasigroup}
A set $C\subset {\Sigma}^n$ is an $(n,2)_4$ MDS code if and only if
$ C=\{(\bar x,g(\bar x))\,|\,\bar x\in {\Sigma}^{n-1}\} $
for some $(n-1)$-qua\-si\-group $g$.
\epro
The following statements are also well known and easy to prove.
\bpro
{\rm (a)} The superposition $g(\bar x^{\laaa i\raaa}\lbb f(\bar y)\rbb)$ of an $n$-qua\-si\-group $g$ and
          an $m$-qua\-si\-group $f$ is an $(n+m-1)$-qua\-si\-group.\\
{\rm (b)} If $g$ is an $n$-qua\-si\-group and $i\in [n]$,
          then its inversion in the $i$th position $f^{\langle i\rangle}$ is an $n$-qua\-si\-group too.\\
{\rm (c)} If $g$ is an $n$-qua\-si\-group and $a\in {\Sigma}$,
          then the set
          \[M_{a}\triangleq \{\bar x\in {\Sigma}^n \,|\, g(\bar x)=a\}\]
          is an MDS code.\\
{\rm (d)} A $1$-qua\-si\-group $p:{\Sigma} \to {\Sigma}$ is a permutation of ${\Sigma}$.
\epro

A set $S\subset {\Sigma}^n$ is called a {\em double-code\/} if each line
of ${\Sigma}^n$ contains zero or two elements from $S$. A double-code
$S\subset {\Sigma}^n$ is called {\em double-MDS-code\/} if each line of
${\Sigma}^n$ contains exactly two elements from $S$. If a double-code is
a subset of some double-MDS-code, then we call it \em complementable\em.
If a double-code is complementable, nonempty, and cannot be
split into more than one nonempty double-codes, then we call
it \em prime\em.

{\em Remark}. The union of two disjoint
 $(n,
 2)_4$ MDS codes is always
a double-MDS-code. The converse statement does not hold for $n\geq 3$
(see e.\,g. {\rm \cite{KroPot:nonsplittable}}).

\begin{figure}[hb]
\renewcommand\arraystretch{0.8}
\mbox{}\hfill {\rm a)}
\begin{tabular}[c]{|@{\hspace{0.6ex}}c@{\hspace{0.6ex}}|@{\hspace{0.6ex}}c@{\hspace{0.6ex}}|@{\hspace{0.6ex}}c@{\hspace{0.6ex}}|@{\hspace{0.6ex}}c@{\hspace{0.6ex}}|}
\hline \phantom{\raisebox{-0.2ex}{$\circ$}}&\phantom{\raisebox{-0.2ex}{$\circ$}}&\phantom{\raisebox{-0.2ex}{$\circ$}}&\phantom{\raisebox{-0.2ex}{$\circ$}}\\
\hline \phantom{\raisebox{-0.2ex}{$\circ$}}&\phantom{\raisebox{-0.2ex}{$\circ$}}&\phantom{\raisebox{-0.2ex}{$\circ$}}&\phantom{\raisebox{-0.2ex}{$\circ$}}\\
\hline \phantom{\raisebox{-0.2ex}{$\circ$}}&\phantom{\raisebox{-0.2ex}{$\circ$}}&\phantom{\raisebox{-0.2ex}{$\circ$}}&\phantom{\raisebox{-0.2ex}{$\circ$}}\\
\hline \phantom{\raisebox{-0.2ex}{$\circ$}}&\phantom{\raisebox{-0.2ex}{$\circ$}}&\phantom{\raisebox{-0.2ex}{$\circ$}}&\phantom{\raisebox{-0.2ex}{$\circ$}}\\
\hline
\end{tabular}
\hfill {\rm b)}
\begin{tabular}[c]{|@{\hspace{0.6ex}}c@{\hspace{0.6ex}}|@{\hspace{0.6ex}}c@{\hspace{0.6ex}}|@{\hspace{0.6ex}}c@{\hspace{0.6ex}}|@{\hspace{0.6ex}}c@{\hspace{0.6ex}}|}
\hline\phantom{\raisebox{-0.2ex}{$\circ$}}&\phantom{\raisebox{-0.2ex}{$\circ$}}&\phantom{\raisebox{-0.2ex}{$\circ$}}&\phantom{\raisebox{-0.2ex}{$\circ$}}\\
\hline\phantom{\raisebox{-0.2ex}{$\circ$}}&\phantom{\raisebox{-0.2ex}{$\circ$}}&\phantom{\raisebox{-0.2ex}{$\circ$}}&\phantom{\raisebox{-0.2ex}{$\circ$}}\\
\hline\raisebox{-0.2ex}{$\bullet$}&\raisebox{-0.2ex}{$\bullet$}&\raisebox{-0.2ex}{$\circ$}&\raisebox{-0.2ex}{$\circ$}\\
\hline\raisebox{-0.2ex}{$\bullet$}&\raisebox{-0.2ex}{$\bullet$}&\raisebox{-0.2ex}{$\circ$}&\raisebox{-0.2ex}{$\circ$}\\
\hline
\end{tabular}
\hfill {\rm c)}
\begin{tabular}[c]{|@{\hspace{0.6ex}}c@{\hspace{0.6ex}}|@{\hspace{0.6ex}}c@{\hspace{0.6ex}}|@{\hspace{0.6ex}}c@{\hspace{0.6ex}}|@{\hspace{0.6ex}}c@{\hspace{0.6ex}}|}
\hline \raisebox{-0.2ex}{$\circ$}&\raisebox{-0.2ex}{$\circ$}&\raisebox{-0.2ex}{$\bullet$}&\raisebox{-0.2ex}{$\bullet$}\\
\hline \raisebox{-0.2ex}{$\circ$}&\raisebox{-0.2ex}{$\circ$}&\raisebox{-0.2ex}{$\bullet$}&\raisebox{-0.2ex}{$\bullet$}\\
\hline \raisebox{-0.2ex}{$\bullet$}&\raisebox{-0.2ex}{$\bullet$}&\raisebox{-0.2ex}{$\circ$}&\raisebox{-0.2ex}{$\circ$}\\
\hline \raisebox{-0.2ex}{$\bullet$}&\raisebox{-0.2ex}{$\bullet$}&\raisebox{-0.2ex}{$\circ$}&\raisebox{-0.2ex}{$\circ$}\\
\hline
\end{tabular}
\hfill {\rm d)}
\begin{tabular}[c]{|@{\hspace{0.6ex}}c@{\hspace{0.6ex}}|@{\hspace{0.6ex}}c@{\hspace{0.6ex}}|@{\hspace{0.6ex}}c@{\hspace{0.6ex}}|@{\hspace{0.6ex}}c@{\hspace{0.6ex}}|}
\hline \raisebox{-0.2ex}{$\circ$}&\raisebox{-0.2ex}{$\circ$}&\raisebox{-0.2ex}{$\bullet$}&\raisebox{-0.2ex}{$\bullet$}\\
\hline \raisebox{-0.2ex}{$\circ$}&\raisebox{-0.2ex}{$\bullet$}&\raisebox{-0.2ex}{$\circ$}&\raisebox{-0.2ex}{$\bullet$}\\
\hline \raisebox{-0.2ex}{$\bullet$}&\raisebox{-0.2ex}{$\circ$}&\raisebox{-0.2ex}{$\bullet$}&\raisebox{-0.2ex}{$\circ$}\\
\hline \raisebox{-0.2ex}{$\bullet$}&\raisebox{-0.2ex}{$\bullet$}&\raisebox{-0.2ex}{$\circ$}&\raisebox{-0.2ex}{$\circ$}\\
\hline
\end{tabular}
\hfill {\rm e)}
\begin{tabular}[c]{|@{\hspace{0.6ex}}c@{\hspace{0.6ex}}|@{\hspace{0.6ex}}c@{\hspace{0.6ex}}|@{\hspace{0.6ex}}c@{\hspace{0.6ex}}|@{\hspace{0.6ex}}c@{\hspace{0.6ex}}|}
\hline & & & \\
\hline \raisebox{-0.2ex}{$\circ  $}&\raisebox{-0.2ex}{$\bullet$}&\raisebox{-0.2ex}{$\bullet$}&\raisebox{-0.2ex}{$\circ$} \\
\hline \raisebox{-0.2ex}{$\bullet$}&\raisebox{-0.2ex}{$\circ  $}&\raisebox{-0.2ex}{$\bullet$}&\raisebox{-0.2ex}{$\circ$} \\
\hline \raisebox{-0.2ex}{$\bullet$}&\raisebox{-0.2ex}{$\bullet$}&\raisebox{-0.2ex}{$\circ  $}&\raisebox{-0.2ex}{$\circ$} \\
\hline
\end{tabular}\hfill\mbox{}
\caption{$\bullet$ -- the elements of double-codes in ${\Sigma}^2$;\ \ \
$\circ$ -- the result of the operation $\backslash_1$.}

\end{figure}

\bexam\label{ex1} {\rm Figure 1 shows all double-codes in ${\Sigma}^2$ up
to permutations of rows and columns. The double-codes a)-d) are
complementable and e) is not. The double-codes c) and d) are
double-MDS-codes. The double-codes b) and d) are prime. }
\eexam
\section{Preliminary statements}\label{sect:pre}
\mbox{}

\bpro \label{p1}
{\rm (a)} If $S\subset {\Sigma}^n$ is a double-MDS-code, then
its supplement ${\Sigma}^n \setminus S$ is a double-MDS-code. \\
{\rm (b)} A double-code $S\subset {\Sigma}^n$
is a double-MDS-code if and
only if $|S|=|{\Sigma}^n|/2=2^{2n-1}$.
\epro
\proofr (a) follows from the
definition of a double-MDS-code.
(b) is obvious if we consider the partition of ${\Sigma}^n$ into
$i$-lines where $i\in [n]$ is fixed. \proofend

For arbitrary subset $S\subseteq {\Sigma}^n$ we define the {\it adjacency
graph} $G(S)$ with  vertex set $S$, where two vertices are adjacent if
and only if they differ in exactly one coordinate.

The following proposition gives a natural treatment of a complementable double-code
in terms of connected components of the adjacency graph of a double-code
that includes the given double-code.
\bpro\label{p1-3}
Let $S$ be a complementable double-code and $S_0$
be an arbitrary subset of $S$; then \\
{\rm (a)} $S_0$ is a double-code if and only if
$G(S_0)$ is a union of connected components of $G(S)$;\\
{\rm (b)} $S_0$ is a prime double-code if and only if $G(S_0)$ is
a connected component of $G(S)$.
\epro
\proofr
The graph $G(S)$ has an edge between $S_0$ and $S\backslash S_0$
if and only if there is a line that has nonempty intersections
with both $S_0$ and $S\backslash S_0$.
Now, (a) follows from the definitions of double-codes and connected components of a graph.
(b) can be easily derived from (a).
%
%
%
\proofend
\bcor\label{c1-3}
Assume that prime double-codes $C$ and $C'$ are included in the same
complementable double-code. Then $C=C'$ or $C\cap C'=\emptyset$.
\ecor

The following simple proposition
will be used in Sections~\ref{sect:MDS} and~\ref{sect:quasigroup}.
\bpro\label{p0} Let $S$ be a double-MDS-code and let $\gamma$ be the
number of prime double-codes included in $S$. {\rm (a)} If $G(S)$ is a
bipartite graph, then $S$ includes exactly $2^{\gamma}$
different MDS codes. {\rm (b)} Otherwise, $S$ does not include an MDS code.
\epro
\proofr
(a) By Proposition~\ref{p1-3}(b), $\gamma$ is the number of connected components
in $G(S)$.
A part of the bipartite graph $G(S)$ is an MDS code by the definition.
So,  the number of the MDS codes that $S$ includes equals the number of the ways of choosing
a  part of the bipartite graph $G(S)$, i.\,e., $2^{\gamma}$.

(b) Assume that a double-MDS-code $S$ includes an MDS code $C$.
Then, by the definition, $S\backslash C$ also is an MDS code.
So, the graphs $G(C)$ and $G(S\backslash C)$ do not contain edges,
and hence the graph $G(S)$ is bipartite
.
\proofend

Let $S\subset {\Sigma}^n$ and $i\in [n]$; then we denote
\[
\displaystyle\backslash_i S \triangleq {\mathcal E}_i(S)\backslash S
\]
(see Fig.~1 for example).

\bpro\label{abcdefghijk}
Let $S,S'\subset {\Sigma}^n$ be double-codes and $i,i'\in [n]$. Then\\
{\rm (a)} $S\cap\backslash_i S=\emptyset$;\\
{\rm (b)} $\backslash_i \backslash_i S=S$;\\
{\rm (c)} $|S|=|\backslash_i S|$;\\
{\rm (d)} $S\subseteq S'$ if and only if $\backslash_i S\subseteq \backslash_i S'$;\\
{\rm (e)} if $S \cup S'$ is a double-code, then
          $\backslash_i (S \cup S')=\backslash_i S \cup \backslash_i S'$;\\
{\rm (f)} $S$ is a double-MDS-code if and only if $\backslash_i S$ is a double-MDS-code;\\
\phantom{\rm (f)} $S$ is a double-MDS-code if and only if $\backslash_i S={\Sigma}^n\backslash S$;\\
{\rm (g)} $S$ is complementable if and only if $\backslash_i S$  is a complementable double-code;\\
{\rm (h)} $S$ is prime if and only if $\backslash_i S$ is a prime double-code;\\
{\rm (i)} if $S$  is prime, then either
   $\backslash_i S=\backslash_{i'}S$ or $\backslash_i S\cap \backslash_{i'}S=\emptyset$;\\
{\rm (j)} if $S$ is complementable, then
          $\backslash_i \backslash_{i'} S=\backslash_{i'}\backslash_i S$;\\
{\rm (k)} $S$ is a double-MDS-code if and only if $|S|>0$ and $\backslash_j S
= \backslash_{j'}S$ for each $j,j'\in [n]$.
\epro

\proofr (a) is clear. The set ${\mathcal E}_i(S)={\mathcal
E}_i(\backslash_i S)=S\cup \backslash_i S$ can be partitioned into
$i$-lines. Each line of the partition has two elements from $S$
and the other two from $\backslash_i S$. Now (b) and (c) are also obvious.

(d) Suppose, $S\subseteq S'$. Then ${\mathcal E}_i(S)\subseteq {\mathcal
E}_i(S')$. Each line ${\mathcal E}_i(\bar x)$,
where $\bar x\in {\mathcal E}_i(S)$,
contains two elements from $S$ and the other two from
$\backslash_{i} S$. They also are elements of $S'$ and
$\backslash_{i}S'$ respectively. So, each element from
$\backslash_{i} S$ is in $\backslash_{i} S'$. The converse
statement is proved in the same way.

(e) It is easy to see that each $i$-line $\mathcal E$  \\
- either has the same intersection with both $S$ and $S'$\\
- or is disjoint with $S$ or $S'$.\\
In any case,
${\mathcal E}\cap\backslash_i (S \cup S')
={\mathcal E}\cap(\backslash_i S \cup \backslash_i S')$.
Since ${\Sigma}^n$ is the union of $i$-lines, the statement is proved.

(f) follows from (a), (c), and Proposition~\ref{p1}(a,b).

(g) Assume the double-code $S$ is complementable. First we will show that $\backslash_i S$ is a double-code.
Let ${\mathcal E}_j(\bar x)$ be an arbitrary line, where
$j\in [n]$ and $\bar x=(x_1,\ldots ,x_n)\in {\Sigma}^n$. If $j=i$, then
$|{\mathcal E}_j(\bar x)\cap S|=|{\mathcal E}_j(\bar x)\cap\backslash_i
S|\in\{0,2\}$. Let $j\neq i$. It is clear that ${\mathcal F}_{i,j;\bar
x}S$ is a double-code in ${\Sigma}^2$ and ${\mathcal F}_{i,j;\bar
x}\backslash_i S= \backslash_1 {\mathcal F}_{i,j;\bar x}S$.
Furthermore, the fact that $S$ is complementable implies that ${\mathcal
F}_{i,j;\bar x}S$ is complementable too. It is easy to check (see Fig. 1(a-d))
that $\backslash_1 {\mathcal F}_{i,j;\bar x}S$ is a
double-code. Consequently, $|{\mathcal E}_j(\bar x)\cap\backslash_i
S|=|{\mathcal E}_2(x_i,x_j)\cap\backslash_1 {\mathcal F}_{i,j;\bar
x}S|\in\{0,2\}$ and $\backslash_i S$ is a double-code by the
definition.

Since $S$ is a complementable double-code, there is a double-MDS-code $S''\supseteq S$.
By (f), the set $\backslash_i S''$ is a double-MDS-code.
By (d), we have $\backslash_i S \subseteq \backslash_i S''$.
Consequently, the double-code $\backslash_i S$ is complementable.

Similarly, if $\backslash_i S$ is a complementable double-code,
then $S$ is.

(h) By (g), we may assume that $S$ and $\backslash_i S$ are
complementable double-codes. Let $S$ be non prime, i.\,e.,
$S=S_1\cup S_2$, where $S_1$ and $S_2$ are disjoint nonempty
double-codes. Double-codes $S_1$ and $S_2$ are complementable by the definition;
$\backslash_i S_1$ and $\backslash_i S_2$ are also complementable double-codes by (g).
The sets ${\mathcal E}_i(S_1)$ and ${\mathcal E}_i(S_2)$ are disjoint.
Therefore, $\backslash_i S_1$ and $\backslash_i S_2$ are disjoint
and the double-code $\backslash_i S=\backslash_i
S_1\cup\backslash_i S_2$ is not prime.
This proves that if $\backslash_i S$ is prime, then $S$ is prime.
Similarly, the converse also holds.

(i) Let $S\subseteq S''$, where $S''$ is a double-MDS-code. It
follows from (d) and (f) that $\backslash_i S\subseteq \backslash_i
S''={\Sigma}^n\backslash S''$. On the other hand, $\backslash_{i'}
S\subseteq \backslash_{i'} S''={\Sigma}^n\backslash S''$.
By (h), the sets $\backslash_i S$ and $\backslash_{i'} S$ are prime double-codes;
by Corollary~\ref{c1-3}, they are either coincident or disjoint.

(j) It is enough to check that
for each $\bar x\in {\Sigma}^n$ it holds
$\backslash_i \backslash_{i'}S_{i,i';\bar x}=\backslash_{i'}\backslash_i S_{i,i';\bar x}$,
where $S_{i,i';\bar x}\triangleq S\cap\{\bar x^{\laaa i,i'\raaa}\lbbb b,c\rbbb|b,c\in {\Sigma}\}$.
Equivalently, $\backslash_1 \backslash_2{\mathcal F}_{i,i';\bar
x}S=\backslash_2\backslash_1{\mathcal F}_{i,i';\bar x}S$ for all $\bar
x\in {\Sigma}^n$. The last can be checked directly, taking into account
that ${\mathcal F}_{i,i';\bar x}S$ is a complementable double-code
(see Fig.\,1a-d).

(k) We first note that the condition $\backslash_j S = \backslash_{j'}S$
for all $j,j'\in [n]$ is equivalent to the condition
$
\backslash_{1}{\mathcal F}_{i,j;\bar x}S=\backslash_{2}{\mathcal F}_{i,j;\bar x}S$
for all different $i,j\in [n]$ and $\bar x\in {\Sigma}^n$.
Since ${\mathcal F}_{i,j;\bar x}S$ is a double-code,
it is straightforward (see Fig.\,1) that the last condition
is equivalent to
\begin{equation}\label{eq:12F}
{\mathcal F}_{i,j;\bar x}S=\emptyset \qquad\mbox{or}\qquad
{\mathcal F}_{i,j;\bar x}S\mbox{ is a double-MDS-code.}
\end{equation}

\emph{Only if:} If $S$ is a double-MDS-code,
then
(\ref{eq:12F}) holds automatically.

\emph{If:} Suppose, by contradiction, that $S$ is not a double-MDS-code.
Then there exist $\bar v=(v_1,\ldots ,v_n)$ and
$\bar z=(z_1,\ldots ,z_n)\in {\Sigma}^n$ such that ${\mathcal E}_1(\bar v)\cap
S\neq\emptyset$ and ${\mathcal E}_1(\bar z)\cap S=\emptyset$. Consider
the sequence $\bar v=\bar v^0,\bar v^1,\ldots ,\bar v^n=\bar z$, where
$\bar v_j=(z_1,\ldots ,z_j,v_{j+1},\ldots ,v_n)$. Note that $\bar v^{j-1}$
and $\bar v^{j}$ coincide in all positions may be except the $j$th one.
There exists $j\in [n]$ such that
${\mathcal E}_1(\bar v^{j-1})\cap S\neq\emptyset$
and
${\mathcal E}_1(\bar v^j)\cap S=\emptyset$.
Then $|{\mathcal F}_{1,j;\bar v^j}S|$
contradicts to (\ref{eq:12F}) with $i=1$ and $\bar x=\bar v^j$.
\proofend
Let $S$ be a double-MDS-code in ${\Sigma}^n$ and let $R\subseteq S$ be a
prime double-code. We say that  $i$ and $i'$ from
$[n]$ are {\it equivalent} (and write $i\stackrel{\scriptscriptstyle S}\sim i'$) if
$\backslash_{i}R=\backslash_{i'}R$.
In Corollary~\ref{p4+} below we will show
that the equivalence ${\stackrel{\scriptscriptstyle S}{\sim}}$ does not depend on the choice of $R$.
Let $K_1=\{i_{1,1},i_{1,2},\ldots ,i_{1,n_1}\},
     K_2=\{i_{2,1},i_{2,2},\ldots ,i_{2,n_2}\},\ldots,
     K_k=\{i_{k,1},i_{k,2},\ldots ,i_{k,n_k}\}$
be the equivalence classes of ${\stackrel{\scriptscriptstyle S}{\sim}}$.
By Proposition~\ref{abcdefghijk}(k), we have the following:
%
\bpro\label{p3}
The double-MDS-code $S$ is prime if and only if $k=1$.
\epro

Denote by ${\{0,1\}^n_{even}}$
the set of the even-weight (i.\,e., with even number of ones) elements of $\{0,1\}^n$.
Denote by $\bar e_j$ the word in $\{0,1\}^n$ with the only $1$ in the $j$th position.
Let the sets $R_{\bar y}$, where $\bar y\in \{0,1\}^n$,
be inductively defined by the equalities
\[
R_{\bar 0} \triangleq R
\qquad\mbox{and}\qquad
R_{\bar y\oplus \bar e_j} \triangleq \backslash_{j}R_{\bar y}.
\]
As follows from Proposition~\ref{abcdefghijk}(b,j),
the sets $R_{\bar y}$ are well defined.
The next proposition is a corollary of Proposition~\ref{abcdefghijk}.

\bpro\label{p4}\mbox{}\\
{\rm (a)} For each $\bar y\in \{0,1\}^n$
    the set $R_{\bar y}$ is a prime double-code.\\
{\rm (b)} For each $\bar y\in\{0,1\}^n$ the equality
$\backslash_{i'} R_{\bar y}=\backslash_{i''} R_{\bar y}$
holds if and only if $i' \stackrel{\scriptscriptstyle S}\sim i''$.\\
{\rm (c)} For each $\bar y, \bar z\in \{0,1\}^n$
either $R_{\bar y}=R_{\bar z}$ or $R_{\bar y}\cap R_{\bar z}=\emptyset$. \\
{\rm (d)} $S=\bigcup_{\bar y\in {\{0,1\}^n_{even}}} R_{\bar y}$. \epro
\proofr (a) follows from Proposition~\ref{abcdefghijk}(h).

(b) Let $\bar y=\bar e_{j_1}\oplus\ldots\oplus \bar e_{j_w}$. Then $R_{\bar
y}=\backslash_{j_w}{\cdot}{\cdot}{\cdot}\backslash_{j_1} R$. By Proposition
\ref{abcdefghijk}(b,j), we have
\begin{eqnarray*}\hspace{-3ex}
\backslash_{i'} R_{\bar y}=\backslash_{i''} R_{\bar y} &
\Longleftrightarrow & \backslash_{i'}
\backslash_{j_w}{\cdot}{\cdot}{\cdot}\backslash_{j_1} R=
\backslash_{i''}\backslash_{j_w}{\cdot}{\cdot}{\cdot}\backslash_{j_1} R
\\ & \Longleftrightarrow &
\backslash_{j_w}{\cdot}{\cdot}{\cdot}\backslash_{j_1} \backslash_{i'} R=
\backslash_{j_w}{\cdot}{\cdot}{\cdot}\backslash_{j_1} \backslash_{i''}R
\Longleftrightarrow \backslash_{i'} R=\backslash_{i''}R
\Longleftrightarrow i' \stackrel{\scriptscriptstyle S}\sim i''.
\end{eqnarray*}

(c,d) Let again $R_{\bar y}=\backslash_{j_w}{\cdot}{\cdot}{\cdot}\backslash_{j_1} R$.
We have $R\subseteq S$. It follows by induction from Proposition
\ref{abcdefghijk}(d,f) that $R_{\bar y}\subseteq S$ if $\bar
y\in\{0,1\}^n_{even}$ and $R_{\bar y}\subseteq {\Sigma}^n\backslash S$
otherwise.
So, (c) holds by Proposition~\ref{p1}(a) and Corollary~\ref{c1-3}.
Further, the union $\bigcup_{\bar y\in {\{0,1\}^n_{even}}}
R_{\bar y}$ is a subset of $S$ and, by the definition, is a complementable double-code.
On the other hand, by Proposition~\ref{abcdefghijk}(e), the set
\[ \backslash_{i}\bigcup_{\bar y\in {\{0,1\}^n_{even}}} R_{\bar y}
= \bigcup_{\bar y\in {\{0,1\}^n_{even}}} \backslash_{i}R_{\bar y}
= \bigcup_{\bar y\in {\{0,1\}^n_{even}}} R_{\bar y\oplus \bar e_i}
= \bigcup_{\bar y\in {\{0,1\}^n_{odd}}} R_{\bar y} \]
does not depend on $i$ and, by Proposition~\ref{abcdefghijk}(k), the set
$\bigcup_{\bar y\in {\{0,1\}^n_{even}}} R_{\bar y}$
is a double-MDS-code. Therefore, it coincides with $S$. \proofend

\bcor\label{p4+} The equivalence ${\stackrel{\scriptscriptstyle S}{\sim}}$ does not depend on the
choice of the prime double-code $R\subseteq S$. \ecor

\proofr It follows from Proposition~\ref{p4}(d) and Corollary~\ref{c1-3} that
for each prime double-code $R'\subseteq S$ there exists
${\bar y\in {\{0,1\}^n_{even}}}$ such that $R'=R_{\bar y}$.
Proposition~\ref{p4}(b) completes the proof. \proofend

\bcor \label{p4++}
Let $S,S',S''\subset {\Sigma}^n$ be double-MDS-codes
and $S_0$ be a double-code.\\
{\rm(a)} If $S_0 \subseteq S\cap S'$, then $S =S'$.
\\
{\rm(b)} If $S_0 \subseteq S\backslash S'' $, then $S={\Sigma}^n\backslash S''$.
\ecor
\proofr
Let $R\subseteq S_0$ be a
prime double-code. Then $R\subseteq S$ and $R\subseteq S'$. By
Proposition~\ref{p4}(d), we have $S=\bigcup_{\bar y\in
{\{0,1\}^n_{even}}} R_{\bar y}$ and $S'=\bigcup_{\bar y\in
{\{0,1\}^n_{even}}} R_{\bar y}$.
So, (a) is proved;
taking $S' ={\Sigma}^n\backslash S''$, we get (b). \proofend

Let $\sigma=\chi_S:{\Sigma}^n\to \{0,1\}$ be the characteristic function
of the double-MDS-code $S$ and let for each $j\in [k]$
$$ \sigma_j(y_1,\ldots ,y_{n_j})
   \triangleq \sigma(\bar 0^{\laaa i_{j,1},\ldots ,i_{j,n_j}\raaa}\lbbb y_1,\ldots ,y_{n_j}\rbbb) $$
be its subfunction with the set of arguments that corresponds to the class
$K_j$.

\bpro \label{p5} For any nonequivalent $i',i''\in [n]$,
any $\bar x\in {\Sigma}^n$, and any $a',a''\in {\Sigma}$
it holds
$$ \sigma(\bar x) \oplus \sigma(\bar x^{\laaa i'\raaa}\lbb a'\rbb)
\oplus \sigma(\bar x^{\laaa i''\raaa}\lbb a''\rbb)
\oplus \sigma(\bar x^{\laaa i',i''\raaa}\lbbb a',a''\rbbb)=0.
$$
\epro

\proofr
Let $\bar y\in {\{0,1\}^n_{even}}$ be such that
$R_{\bar y}\cap {\mathcal E}_{i'}(\bar x)\neq\emptyset$
(by Proposition~\ref{p4}(d), such $\bar y$ exists).
Let us consider the double-MDS-code $S^2={\mathcal F}_{i',i'';\bar x}S\subset {\Sigma}^2$
and the prime double-code
$R^2_{\bar y}={\mathcal F}_{i',i'';\bar x}R_{\bar y}\subset S^2$.
Since, by Proposition~\ref{p4}(b),
$\backslash_{i'}R_{\bar y}\neq\backslash_{i''}R_{\bar y}$,
we find by Proposition~\ref{abcdefghijk}(i) that
$\backslash_{i'}R_{\bar y}\cap\backslash_{i''}R_{\bar y}=\emptyset$
and, consequently,
$\backslash_{1}R_{\bar y}^2\cap\backslash_{2}R_{\bar y}^2=\emptyset$.
Therefore, $R_{\bar y}^2$ corresponds to the case b) of Fig.~1
and $S^2$ corresponds to the case c),
up to permutations of the rows and the columns.
For this case, it is easy to check that
\[
\chi_{S^2}(b',b'') \oplus
\chi_{S^2}(a',b'')  \oplus \chi_{S^2}(b',a'') \oplus
\chi_{S^2}(a',a'')=0\ \ \ \forall b',b'',a',a''\in {\Sigma}.
\]
The statement follows from the obvious identity

$\qquad\sigma(\bar x^{\laaa i',i''\raaa}\lbbb c',c''\rbbb) =
\chi_{S^2}(c',c'')\ \ \ \forall c',c''\in {\Sigma}$.
\hfill\proofend

In the proofs of the following statements we will use
the results and notation
of the Appendix on functions with separable arguments.

\bpro \label{p7}
For each $\bar x$ from ${\Sigma}^n$ it holds
\begin{equation}\label{eq:si1k}
\sigma(\bar x)\equiv\bigoplus_{j=1}^k
\sigma_j(x_{i_{j,1}},x_{i_{j,2}},\ldots ,x_{i_{j,n_j}}) \oplus
\sigma_0,\quad\mbox{where }\sigma_0=(k-1)\sigma(\bar 0).
\end{equation}
\epro

\proofr
By the criterion of Lemma~\ref{l1} of the Appendix,
Proposition~\ref{p5} means that $\sigma$ has
$\{K_1,\ldots,K_k\}$-separable arguments, i.\,e.,
$$
\sigma(\bar x) \equiv \bigoplus_{j=1}^k
f_j(x_{i_{j,1}},x_{i_{j,2}},\ldots ,x_{i_{j,n_j}})
$$
for some functions $f_j:\Sigma^{n_j}\to\{0,1\}$. Then,
\begin{equation}
\label{p7e1}
\sigma(\bar x)\oplus \sigma(\bar 0) \equiv \bigoplus_{j=1}^k
\left(f_j(x_{i_{j,1}},x_{i_{j,2}},\ldots ,x_{i_{j,n_j}})\oplus
f_j(\bar 0)\right).
\end{equation}
Setting $x_i=0$ for all $i\notin K_j$, we have
\begin{equation}\label{p7e2}
\sigma_j(x_{i_{j,1}},x_{i_{j,2}},\ldots ,x_{i_{j,n_j}})\oplus
\sigma(\bar 0) \equiv
f_j(x_{i_{j,1}},x_{i_{j,2}},\ldots ,x_{i_{j,n_j}})\oplus f_j(\bar 0).
\end{equation}
Substituting (\ref{p7e2}) to (\ref{p7e1}) proves the statement.
\proofend
\bpro \label{p8}
For each $j\in [k]$ the function $\sigma_j$
is the characteristic function of a prime double-MDS-code.
\epro

\proofr
For each $j\in [k]$
the function $\sigma_j$ is a subfunction of $\sigma$
and, consequently, the characteristic function of some double-MDS-code $S_j$.
It remains to prove that $S_j$ is prime.
The idea is to show that the $\stackrel{\scriptscriptstyle S}\sim$ equivalence
of the indexes from $K_j$ yields
the $\stackrel{\scriptscriptstyle S_j}\sim$ equivalence of the indexes from $[n_j]$.

We first observe the following straightforward fact:\\
(*) if $R_1\subset \Sigma^{n_1}$, \ldots, $R_k\subset \Sigma^{n_k}$
are double codes,
\begin{equation}\label{eq:R1k}
 R\triangleq R_1\times\ldots\times R_k,
\end{equation}
 and $i\in[n_1]$, then $R$ is a double-code and
\( \backslash_i R=(\backslash_i R_1)\times R_2\times\ldots\times R_k. \)

Without loss of generality assume that
$K_1=\{1,\ldots,n_1\}$, $K_2=\{n_1+1,\ldots,n_1+n_2\}$, and so on.
Then, from (\ref{eq:si1k}) we derive that
$S\supseteq S_1\times\ldots\times S_k$ or
$\Sigma^n\setminus S\supseteq S_1\times\ldots\times S_k$.
For each $j\in [k]$ we choose a prime double-code $R_j\subseteq S_j$.
Then, the double-code $R$ defined as in (\ref{eq:R1k})
is included in $S$ or $\Sigma^n\setminus S$,
and for any equivalent (in the sense of $\stackrel{\scriptscriptstyle S}\sim$)
$i$ and $i'$ we have $\backslash_{i} R=\backslash_{i'} R$.
From (*) we derive that for each $i,i'\in K_1$
we have $\backslash_{i} R_1=\backslash_{i'} R_1$,
i.\,e., $i \stackrel{\scriptscriptstyle S_1}\sim i'$.
Thus, by Proposition~\ref{p3}, $S_1$ is a prime double-MDS-code.
The same is true of $S_j$ for every $j\in [k]$.
%
\proofend

\section{Decomposition of double-MDS-codes}
\label{sect:2MDS}
\mbox{}

\btheorem \label{th1}
{\rm (a)} The characteristic function $\chi_S$ of a
double-MDS-code $S$ has a unique representation in the form
\begin{equation}\label{e1}
\chi_S(\bar x)=\bigoplus_{j=1}^k \chi_{S_j}(\tilde x_j)\oplus
\sigma_0 \qquad \mbox{where}
\end{equation}
\begin{itemize}
\item
 $k\in [n]$,
\item
 $\tilde x_1$, \ldots, $\tilde x_k$
 are disjoint collections of variables from $\bar x$,
 $\tilde x_j\triangleq(x_{i_{j,1}},\ldots ,x_{i_{j,n_j}})$,
\item
 for each $j\in [k]$ the set $S_j\subset{\Sigma}^{n_j}$ is a prime double-MDS-code
 and $\bar 0\in S_j$,
\item
 $\sigma_0\in \{0,1\}$.
\end{itemize}
{\rm (b)} $S$ is a union of $2^{k-1}$ equipotent prime double-codes;\\
${\Sigma}^n\backslash S$ is a union of $2^{k-1}$ equipotent prime double-codes. \\
{\rm (c)} If $k\geq 2$ and the adjacency graph $G(S)$ is bipartite,
then for each $j\in [k]$ the graphs $G(S_j)$ and $G({\Sigma}^{n_j}\backslash S_j)$
are also bipartite.
\etheorem
\proofr
(a) is a corollary of
Propositions~\ref{p7} and~\ref{p8}.
The uniqueness of the representation (\ref{e1}) follows from
Lemma~\ref{l4} of the Appendix.

(b) Without loss of generality we assume
that the variables are arranged in such a way that
$\bar x=(\tilde x_1,\ldots ,\tilde x_k)$.
Then,
$$
S=\bigcup_{(\gamma_1,\ldots ,\gamma_k)\in \{0,1\}^k\atop
\gamma_1\oplus \ldots \oplus\gamma_k\oplus\sigma_0=0}
S_1^{\gamma_1}\times\ldots\times S_k^{\gamma_k},
\qquad
{\Sigma}^n\backslash S=\bigcup_{(\gamma_1,\ldots ,\gamma_k)\in \{0,1\}^k\atop
\gamma_1\oplus \ldots \oplus\gamma_k\oplus\sigma_0=1}
S_1^{\gamma_1}\times\ldots\times S_k^{\gamma_k}
$$
where $S_j^{0}\triangleq S_j$ and
$S_j^{1}\triangleq {\Sigma}^{n_j}\backslash S_j$. For all
$(\gamma_1,\ldots ,\gamma_k)\in \{0,1\}^k$
the adjacency graph
$G_{(\gamma_1,\ldots ,\gamma_k)}
\triangleq G(S_1^{\gamma_1}\times\ldots\times S_k^{\gamma_k})
=G(S_1^{\gamma_1})\times\ldots\times G(S_k^{\gamma_k})$
is connected and has the degree $n=n_1+\ldots+n_k$, because the graphs
$G(S_1^{\gamma_1}),\ldots,G(S_k^{\gamma_k})$
are connected and have the degrees $n_1,\ldots,n_k$.
Consequently, 
$G_{(\gamma_1,\ldots ,\gamma_k)}$
is a connected component
of $G(S)$     (if $\gamma_1\oplus \ldots \oplus\gamma_k\oplus\sigma_0=0$)
or $G({\Sigma}^n\backslash S)$ (if $\gamma_1\oplus \ldots \oplus\gamma_k\oplus\sigma_0=1$).
Moreover, the cardinality of $G_{(\gamma_1,\ldots ,\gamma_k)}$ equals
\(
|S_1|\cdot\ldots\cdot |S_k|
\)
and does not depend
on $(\gamma_1,\ldots ,\gamma_k)$. Proposition~\ref{p1-3} completes the
proof of (b).

(c) Let $k\geq 2$.
It is easy to see that fixing the arguments
$\tilde x_1, \ldots ,\tilde x_{j-1},\tilde x_{j+1}, \ldots ,\linebreak[1]\tilde x_k$
we can obtain the function $\chi_{S_j}(\tilde x_{j})$,
as well as $\chi_{{\Sigma}^{n_j}\backslash S_j}(\tilde x_j)$,
in the right part of (\ref{e1}).
Consequently,
$\chi_{S_j}$ and $\chi_{{\Sigma}^{n_j}\backslash S_j}$
are subfunctions of $\chi_S$.
Thus, $G({S_j})$ and $G(\Sigma^{n_j}\backslash S_j)$
are subgraphs of $G(S)$, which proves the statement.
%
%
\proofend

\bcor\label{c1} If a double-MDS-code $S$ is not prime, then
$\chi_S$ is the sum of the characteristic functions of prime
double-MDS-codes of smaller lengths.
Moreover, if $G(S)$ has $K$ connected components,
then the number of the summands is $1+\log_2 K$.
\ecor

\bnote {\rm
(On the general $q$-valued case.)
The results above can be generalized to the arbitrary even size of the alphabet
$\Sigma=\{0,1,\ldots,q-1\}$.
A set $S\subset {\Sigma}^n$ is called a {\em $q/2$-fold MDS code\/}
({\em $q/2$-code\/})
if each line of ${\Sigma}^n$ contains $q/2$ (respectively, $0$ or $q/2$)
elements from $S$.
The concepts of \emph{complementable} and \emph{prime} $q/2$-codes are defined as for
double-codes.
If we replace double-codes and double-MDS-codes by, respectively,
$q/2$-codes and $q/2$-fold MDS codes
in Theorem~\ref{th1} and Corollary~\ref{c1},
then the statements will hold as well.
Indeed, all the proofs, without essential changes, are valid for the $q$-valued case.
It should only be noted that for each even $q$ there are exactly one
non-prime $q/2$-fold MDS code in $\Sigma^2$ (see Fig.~1(c) for the case $q=4$)
and exactly one ``non-MDS'' prime $q/2$-code in $\Sigma^2$ (Fig.~1(b)), up to equivalence.
So, it is easy to check that
all the simple statements on the $q/2$-codes
in $\Sigma^2$ that are used in the proofs
(Propositions~\ref{abcdefghijk}(g,j,k) and~\ref{p5})
are valid for the $q$-valued case.
}
\enote

\section{Decomposition of $(n,2)_4$ MDS codes}\label{sect:MDS}
The next theorem gives a representation of
$(n,
2)_4$ MDS codes, which is based on the decomposition of double-MDS-codes
presented in Theorem~\ref{th1}.

A double-MDS-code $S$ in ${\Sigma}^n$ is
called \emph{linear} if its characteristic function $\chi_S$ can
be represented in the form
\[\chi_S(x_1,\ldots,x_n)=\chi_1(x_1)\oplus\ldots\oplus\chi_n(x_n)\]
for some functions $\chi_1,\ldots,\chi_n:\Sigma\to\Sigma$,
which are, clearly, the characteristic functions of double-MDS-codes in ${\Sigma}^1$.
There is only one double-MDS-codes in ${\Sigma}^1$ up to permutation of the symbols of ${\Sigma}$.
So, there is only one linear double-MDS-code in ${\Sigma}^n$ up to permutations of
the alphabet symbols in each coordinate.

An MDS code $C$ is called \emph{semilinear} if $C\subset S$ for
some linear double-MDS-code $S$.
Since in this case $C$ is a part of the bipartite graph $G(S)$,
it is not difficult to describe the class of semilinear MDS codes.
The number of such codes of length $n$ is
$3^{n}2^{2^{n-1}+1}-2^{n+2}3^{n-1}$ (see e.\,g. \cite{PotKro:asymp}).

\btheorem\label{th2}
Let $S$ be a double-MDS-code in ${\Sigma}^n$ and
$C\subset S$ be an $(n,2)_4$ MDS code. Then
\begin{eqnarray}
 C & = & \{(x_1,\ldots,x_n)\,|\,
(g_1(\tilde x_1),\ldots ,g_k(\tilde x_k))\in B_C\},\label{e2} \\
 C & = & \{(x_1,\ldots,x_n)\,|\,
(\tilde x_j,y_j)\in C_j,\ j=1,\ldots ,k;\ (y_1,\ldots ,y_k)\in B_C\},
\label{e3}
\end{eqnarray}
where $\tilde x_j=(x_{i_{j,1}},\ldots ,x_{i_{j,n_j}})$,
the set $B_C$ is a semilinear $(k,2)_4$ MDS code,
the set $C_j$ is a $(n_j+1,2)_4$ MDS code,
the mapping \,$g_j:{\Sigma}^{n_j}\to {\Sigma}$ is a $n_j$-qua\-si\-group,
$j=1,\ldots ,k$, and $k,n_j,i_{j,s}$ are specified by Theorem~\ref{th1}.
Moreover, all the parameters except $B_C$ depend only on $S$
and do not depend on $C\subset S$.
\etheorem

\proofr
It is easy to see that (\ref{e2}) and (\ref{e3}) are equivalent
if $C_j=\{(\tilde z,g_j(\tilde z))\,|\linebreak[1]\,\tilde z\in {\Sigma}^{n_j}\}$.
So, it is enough to show only (\ref{e2}).

If $k=1$, then the statement is obvious. Assume that $k>1$.
By Theorem~\ref{th1}(c),
the graphs $G(S_j)$ and
$G({\Sigma}^{n_j}\backslash S_j)$
are bipartite ($S_j$ are specified in Theorem~\ref{th1}(a)).
Then, for each $j\in [k]$
we can easily define an $n_j$-qua\-si\-group $g_j$ such that its
set of $1$s and $0$s is $S_j$
(more accurately,
define the set of $0$s of $g_j$ as a part of $G(S_j)$,
and the set of $1$s as the other part;
the set of $2$s as a part of $G({\Sigma}^{n_j}\backslash S_j)$,
and the set of $3$s as the other part);
i.\,e.,
\begin{equation}\label{ef:Qj}
    \chi_{S_j}(\tilde x_j) \equiv \chi_{\{0,1\}}(g_j(\tilde x_j)).
\end{equation}

Let the linear double-MDS-code $D\subset {\Sigma}^k$ be defined by the equality
\begin{equation}\label{eq:D-def}
\chi_D(y_1,\ldots ,y_k)=\chi_{\{0,1\}}(y_1)\oplus\ldots\oplus\chi_{\{0,1\}}(y_k)\oplus\sigma_0.
\end{equation}
Using (\ref{ef:Qj}) and (\ref{eq:D-def}),
we can rewrite (\ref{e1}) in the following way:
$$ S=\{(x_1,\ldots ,x_n)\,|\,(g_1(\tilde x_1),\ldots ,g_k(\tilde x_k))\in D\}.
$$
If $B\subset D$ is an MDS code, then the set
$$ \{(x_1,\ldots ,x_n)\,|\,(g_1(\tilde x_1),\ldots ,g_k(\tilde x_k))\in B\}\subset S
$$
is also an MDS code. The double-code $D$ has $2^{2^{k-1}}$
MDS-code-subsets (all these MDS codes are semilinear).
Then, $2^{2^{k-1}}$ different MDS-code-subsets of $S$ are
represented in the form (\ref{e2}).

On the other hand,
by Theorem~\ref{th1}(b), the set $S$ is the union of ${2^{k-1}}$ prime double-codes.
By Proposition~\ref{p0}(a), there are exactly $2^{2^{k-1}}$
subsets of $S$ that are MDS codes. Therefore, all these MDS codes have the
representation (\ref{e2}) and the code $C$ is one of them (with $B=B_C$).
\proofend

We say that an $(n,2)_4$ MDS code $C$ is {\em decomposable} if
there are $m\in \{2,\ldots ,n-2\}$, an $m$-qua\-si\-group $g'$,
an $(n-m)$-qua\-si\-group $g''$, and a permutation $\sigma:[n]\to [n]$ such that
\[C=\{(x_1,\ldots ,x_n)\in\Sigma^n\,|\,
g'(x_{\sigma(1)},\ldots ,x_{\sigma(m)})=g''(x_{\sigma(m+1)},\ldots ,x_{\sigma(n)})\}.\]
Taking into account Proposition~\ref{MDS<->quasigroup},
we can say that a decomposable MDS code can
be represented as a ``concatenation'' of MDS codes of smaller lengths.

\bcor\label{c2}
{\rm (a)}
If $2<k<n$ or $k=2$, $n_1>1$, $n_2>1$, then the MDS code $C$ is decomposable.\\
{\rm (b)} If $k=n$, then $C$ is a semilinear MDS code.
\ecor

\bexam\label{ex2}
\textrm{ Let $\pi=(01)(23)$ and $S=C\cup C'$, where \linebreak
$C'=\{(\pi(x_{1}),x_{2},\ldots ,x_n)|(x_1,\ldots ,x_n)\in C\}$. Then $S$
is a non-prime double-MDS-code, $k\geq 2$, $n_1=1$. The lengths of
the codes $B_C,C_1,\ldots,C_k$ in {\rm(\ref{e3})} are smaller than $n$
(in this case, $C$ and $C'$ are decomposable)
if and only if $2<k<n$. }
\eexam

We say that two sets $C,C'\subseteq {\Sigma}^n$ are {\it isotopic} if
there are permutations $\pi_1,\ldots ,\pi_n:{\Sigma}\to {\Sigma}$ such that
\[
(x_1,\ldots ,x_n)\in C \iff (\pi_1(x_1),\ldots ,\pi_n(x_n))\in C'.
\]
The following theorem means that we cannot get a ``new'' MDS code
if we combine parts of two disjoint MDS codes $C_1$ and $C_2$,
i.\,e., the resulting code can be obtained as semilinear, or
as isotopic to $C_1$ and $C_2$, or it can be composed from MDS codes
of smaller lengths.
\btheorem\label{th:switch}
Let $C_1$ and $C_2$ are disjoint MDS codes.
Suppose an MDS code $C_{new}$ is a subset of $C_1\cup C_2$.
Then there are only three possibilities:\\
either\ {\rm (1)} $C_{new}$ is isotopic with $C_1$ and $C_2$,\\
or {\rm (2)} $C_{new}$ is decomposable,\\
or {\rm (3)} $C_{new}$ is semilinear.
\etheorem
\proofr
Since the MDS codes $C_1$ and $C_2$ are disjoint, the set
$S \triangleq C_1\cup C_2$ is a double-MDS-code.

Consider the representation (\ref{e2}) for the code $C_{new}\subset S$.
By Corollary~\ref{c2}, it is enough to consider the case $k=2$,
$\{n_1,n_2\}=\{1,n-1\}$.
W.\,l.\,o.\,g. we assume that (\ref{e2}) has the form
\[ 
C = \{(x_1,\ldots,x_n)\,|\,
(g_1(x_1),g_2(x_2,\ldots ,x_{n}))\in B_C\}
\] 
with $C=C_{new}$.
By Theorem~\ref{th2}, this equation 
also holds for any MDS code $C\subset S$.
By Proposition~\ref{MDS<->quasigroup}, we have an equivalent equation
\[
C = \{(x_1,\ldots,x_n)\,|\,
f_C(g_1(x_1))=g_2(x_2,\ldots ,x_{n})\}
\]
for any MDS code $C\subset S$, where $1$-qua\-si\-group $f_C$ is defined by
$B_C=\{(y,f_C(y)\,|\linebreak[1] \,y\in {\Sigma}\}$.
Since  $g_1,f_C:{\Sigma}\to {\Sigma}$ are permutations,
any two MDS codes that are subsets of $S$ are isotopic.
\proofend

\section{Decomposition of $n$-qua\-si\-groups of order $4$}\label{sect:quasigroup}
In this section we derive two conditions guaranteing that an $n$-qua\-si\-group
can be represented as a superposition of $n_j$-qua\-si\-groups with $n_j<n$.
Note that, taking into account the one-to-one correspondence
between $n$-qua\-si\-groups and MDS-codes (Proposition~\ref{MDS<->quasigroup}),
the following two theorems are closely related with Theorem~\ref{th2}.

If $f$ is an $n$-qua\-si\-group and $B \subseteq {\Sigma}$, then we denote
$$M_B(f)\triangleq \{\bar x\in {\Sigma}^n\,|\,g(\bar x)\in B\}. $$
\btheorem\label{th_q1}
Let $g$ be an $n$-qua\-si\-group.
Assume ${\Sigma}=\{a,b,c,d\}$ and $S=M_{\{a,b\}}(g)$.
Then
\begin{eqnarray}
 g & = & g_0(g_1(\tilde x_1),\ldots ,g_k(\tilde x_k)),\label{e_q}
\end{eqnarray}
where $\tilde x_j=(x_{i_{j,1}},\ldots ,x_{i_{j,n_j}})$,
the mappings $g_0$, $g_1$, \ldots, $g_k$ are $k$-, $n_0$-,\ldots, $n_k$-qua\-si\-groups,
and $k,n_j,i_{j,s}$ are specified by Theorem~{\rm\ref{th1}}.
\etheorem
\proofr
If $k=1$, then the statement is obvious. Suppose $k>1$.
As in the proof of Theorem~\ref{th2}, we get
$$ S=\{(x_1,\ldots ,x_n)\,|\,(g_1(\tilde x_1),\ldots ,g_k(\tilde x_k))\in D\}
$$
for some linear double-MDS-code $D\in {\Sigma}^k$.

Let $g_0$ be a $k$-qua\-si\-group such that $M_{\{a,b\}}(g_0) = D$. Then
the mapping
\begin{equation}\label{eq:f=q0()}
 f(\bar x)=g_0(g_1(\tilde x_1),\ldots ,g_k(\tilde x_k))
\end{equation}
is an $n$-qua\-si\-group such that $M_{\{a,b\}}(f)=S$.

We claim that:\\
(*) the number of ways to choose $g_0$ equals $2^{2^k}$;\\
(**) the number of $n$-qua\-si\-groups $f$ such that $M_{\{a,b\}}(f)=S$ equals $2^{2^k}$.\\
By Theorem~\ref{th1}, each of the sets
$M_{\{a,b\}}(f)=S$, \ %
$M_{\{c,d\}}(f)={\Sigma}^n \setminus S$, \ %
$M_{\{a,b\}}(g_0)=D$, \ %
$M_{\{c,d\}}(g_0)={\Sigma}^k \setminus D$
is the union of $2^{k-1}$ different prime double-codes.
The number of ways to choose $g_0$ equals
the number of ways to choose an MDS code $M_{\{a\}}(g_0)\subset D$ multiplied by
the number of ways to choose an MDS code $M_{\{c\}}(g_0)\subset {\Sigma}^k \setminus D$,
i.\,e., $2^{2^{k-1}}\cdot 2^{2^{k-1}} = 2^{2^{k}}$
(see Proposition~\ref{p0}(a)).
The claim (*) is proved.
Similarly, (**) is also true.

So, we conclude that all the $n$-qua\-si\-groups $f$ such that $M_{\{a,b\}}(f)=S$
($g$ is one of them) have the representation (\ref{eq:f=q0()}).
\proofend

We say that an $n$-qua\-si\-group $g$ is {\em reducible}
if it can be represented as
a superposition of $n_j$-qua\-si\-groups where $n_j<n$.
We say that an $n$-qua\-si\-group $g$ is {\em semilinear}
if the corresponding MDS code
$\{(\bar x,q(\bar x))\,|\,\bar x\in \Sigma^n\}$ is semilinear
\bcor\label{cor:qua}
Assume the conditions of Theorem~{\rm\ref{th_q1}} holds.
{\rm (a)} If $1<k<n$,
then the $n$-qua\-si\-group $g$ is reducible.
{\rm (b)} If $k=n$, then the $n$-qua\-si\-group $g$ is semilinear.
\ecor
\proofr
(a) is straightforward. (b) From the description of $g_0$
(in the proof of Theorem~{\rm\ref{th_q1}}) we derive that it is semilinear.
Since, in the case $k=n$, the $1$-qua\-si\-groups $g_j$, $j\in[n]$, are just
permutations, $g$ is also semilinear.
\proofend

The next theorem interpret Corollary~\ref{cor:qua} in terms of the $n$-qua\-si\-group
that is inverse to $g$ in some (say, $n$th) argument.

\btheorem
Let $h$ be an $n$-qua\-si\-group, $\{a,b,c,d\}={\Sigma}$,
and $Q$ be the set of $n$-qua\-si\-groups $f$ such that
$f|_{{\Sigma}^{n-1}\times \{a,b\}} \equiv h|_{{\Sigma}^{n-1}\times \{a,b\}}$.
If $2<|Q|<2^{2^{n-1}}$, then the $n$-qua\-si\-group $h$ is reducible.
If $|Q|=2^{2^{n-1}}$, then $h$ is semilinear.
\etheorem
\proofr
Assume $f\in Q$ and $f^{\langle n\rangle}$ denotes the $n$-qua\-si\-group
that is inverse to $f$ in $n$th argument.
Then,
the $a${\small s} and the $b${\small s} of $f^{\langle n\rangle}$ coincide with
respectively $a${\small s} and the $b${\small s} of $h^{\langle n\rangle}$,
i.\,e.,
$M_{\{a\}}(f^{\langle n\rangle})=M_{\{a\}}(h^{\langle n\rangle})$ and $M_{\{b\}}(f^{\langle n\rangle})=M_{\{b\}}(h^{\langle n\rangle})$.
Let $S \triangleq M_{\{a,b\}}(h^{\langle n\rangle})$.
Therefore, $|Q|$ is the number
of ways to choose the $c${\small s} of $f^{\langle n\rangle}$, i.\,e.,
the number of MDS codes $C \subset {\Sigma}^n \setminus S$.
So, by Theorem~\ref{th1}(b) and Proposition~\ref{p0}(a), we have
$|Q|=2^{2^{k-1}}$.
Since by the condition of the theorem $2<|Q|<2^{2^{n-1}}$ (or $|Q|=2^{2^{n-1}}$),
we have $1<k<n$ (respectively, $k=n$).
Then, by Corollary~\ref{cor:qua},
the $n$-qua\-si\-group $h^{\langle n\rangle}$ is reducible (respectively, semilinear).
It is straightforward that $h^{\langle n\rangle}$ is reducible (semilinear)
if and only if $h$ is.
\proofend

\appendix
\section{On functions with separable arguments}

In the Appendix, we will consider the functions with separable arguments,
i.\,e., the functions that can be represented
as the sum of functions of smaller arity depending on mutually disjoint collections
of arguments of the original function.
We will prove a criterion for a function to have separable arguments and will show
that a function has a unique canonical representation as such the sum.

Let ${\Sigma}$ be an arbitrary set that contains $0$.
Let $n$ be a natural number,
$K=\{K_1,\ldots,K_k\}$, where $\emptyset \neq K_j\subseteq [n]$,
be a partition of the set $[n]$,
and $K_j=\{i_{1,1},\ldots,i_{1,n_j}\}$, where $n_j=|K_j|$,
$j\in [k]$.
Let $(\Gamma,\oplus)$ be an Abelian group.
We say that a function $f:{\Sigma}^n\to \Gamma$ has
{\em $K$-separable arguments} if
\begin{equation}\label{defsep}
f(\bar x)\equiv f_1(\bar x_{K_1})\oplus\ldots\oplus f_k(\bar x_{K_k}),
\end{equation}
where $f_j:{\Sigma}^{n_j}\to\Gamma$, $\bar x\triangleq
(x_1,\ldots,x_n)$, and $\bar x_{K_j}\triangleq
(x_{i_{j,1}},\ldots ,x_{i_{j,n_j}})$.
If $|K|>1$, then we say that $f$ has {\em separable arguments}.
We say that a function $f:{\Sigma}^n\to \Gamma$ has {\em non-separable
arguments} if (\ref{defsep}) implies $|K|=1$.
\blemma \label{l1}
A function  $f:{\Sigma}^n\to \Gamma$ has $K$-separable arguments
if and only if for each $i',i''$
that belong to different elements of $K$,
for each $\bar x\in {\Sigma}^n$ and $a',a''\in {\Sigma}$ it holds
\begin{equation}\label{ea1}
f(\bar x) \ominus f(\bar x^{\laaa i'\raaa}\lbb a'\rbb) \ominus f(\bar
x^{\laaa i''\raaa}\lbb a''\rbb) \oplus f(\bar x^{\laaa i',i''\raaa}\lbbb a',a''\rbbb)=0.
\end{equation}
\elemma

\proofr Assume (\ref{ea1}) holds. Let $P=\{p_1,\ldots ,p_m\}\subseteq
[n]$, $Q=\{q_1,\ldots ,q_r\}\subseteq [n]$, and each
$K_j$ is disjoint with at least one of $P$ and $Q$.
Then, by (\ref{ea1}), for each $\bar x=\{x_1,\ldots ,x_n\}\in {\Sigma}^n$ we have \\
\begin{eqnarray*}
\displaystyle & \bigoplus_{s=1}^m \bigoplus_{t=1}^r &
 \bigg(
{f} \left(\bar
0^{\laaa p_1,\ldots ,p_{s-1},q_1,\ldots ,q_{t-1}\raaa}%
\lbbb x_{p_1},\ldots ,x_{p_{s-1}},x_{q_1}, \ldots ,x_{q_{t - 1}}\rbbb \right)
\\ & &
\ominus
{f} \left(\bar
0^{\laaa p_1,\ldots ,p_{s - 1},q_1, \ldots ,q_{t}\raaa}%
\lbbb x_{p_1}, \ldots ,x_{p_{s - 1}},x_{q_1}, \ldots ,x_{q_{t}}\rbbb \right)
\\ & &
 \ominus\,
{f} \left(\bar
0^{\laaa p_1, \ldots ,p_{s},q_1, \ldots ,q_{t - 1}\raaa}
\lbbb x_{p_1}, \ldots ,x_{p_{s}},x_{q_1}, \ldots ,x_{q_{t - 1}}\rbbb \right)
\\ & &
 \oplus
{f} \left(\bar
0^{\laaa p_1, \ldots ,p_{s},q_1, \ldots ,q_{t}\raaa}
\lbbb x_{p_1}, \ldots ,x_{p_{s}},x_{q_1}, \ldots ,x_{q_{t}}\rbbb \right) \bigg)=0.
\end{eqnarray*}
Collecting similar terms we get
\begin{eqnarray}
{f}(\bar 0) \ominus {f}(\bar
0^{\laaa p_1,\ldots ,p_m\raaa}\lbbb x_{p_1},\ldots ,x_{p_m}\rbbb)
\ominus {f}(\bar
0^{\laaa q_1,\ldots ,q_r\raaa}\lbbb x_{q_1},\ldots ,x_{q_r}\rbbb)\nonumber\\\label{ae2}
\oplus {f}(\bar
0^{\laaa p_1,\ldots ,p_m,q_1,\ldots ,q_r\raaa}%
\lbbb [x_{p_1},\ldots ,x_{p_m},x_{q_1},\ldots ,x_{q_r}\rbbb)
&=&0.
\end{eqnarray}

Without loss of generality we can assume that
$\bar x=(\bar x_{K_1},\bar x_{K_2},\ldots ,\bar x_{K_k})$, i.\,e.,
$(i_{1,1},\ldots ,i_{1,n_1},\linebreak[1] i_{2,1},\linebreak[1] \ldots ,i_{2,n_2},\ldots,i_{k,n_k})
 =(1,\ldots,n)$
(recall that $\bar x_{K_j} =(x_{i_{j,1}},\ldots ,x_{i_{j,n_j}})$).
By (\ref{ae2}), it holds
\begin{eqnarray*}
\bigoplus_{j=2}^k {f}(\bar 0)\ominus{f}(\bar x_{K_1},\ldots ,\bar x_{K_{j-1}},0,0\,\ldots \,0)
\ominus {f}(0\,\ldots \,0,\bar x_{K_j},0\,\ldots \,0) \\
\oplus {f} (\bar x_{K_1},\ldots ,\bar x_{K_{j-1}},\bar x_{K_j},0\,\ldots \,0) & = & 0.
\end{eqnarray*}
Collecting similar terms we get
$$ {f}(\bar x)=\bigoplus_{j=1}^k {f}(0,\ldots ,0,\bar x_{K_j},0,\ldots ,0)
\ominus (k-1){f}(\bar 0) $$
and the function $f$ has $K$-separable arguments by the definition.

The inverse statement is straightforward. \proofend

\blemma\label{l2} If a function $f:{\Sigma}^n\to\Gamma$ has
$K$-separable arguments, then the functions $f_1,\ldots,f_k$ such
that
$$ f(\bar x)\equiv\bigoplus_{i=1}^k f_i(\bar x_{K_i}) $$
are uniquely defined up to constant summand.
\elemma
\proofr Let
$$ f(\bar x)\equiv\bigoplus_{i=1}^k f_i(\bar x_{K_i})
   \equiv\bigoplus_{i=1}^k g_i(\bar x_{K_i}). $$
Let $j$ be fixed. Setting $\bar x_{K_i}=(0,\ldots ,0)$ for $i\neq j$, we get
$$ g_j(\bar x_{K_j}) \ominus f_j(\bar x_{K_j})\equiv\bigoplus_{i\neq j}
f_i(0,\ldots ,0)\ominus\bigoplus_{i\neq j} g_i(0,\ldots ,0)\equiv
const \in \Gamma.$$
\proofend


If $K=\{K_i\}_{i=1}^k=\{K_1,\ldots,K_k\}$ and $L=\{L_j\}_{j=1}^l$
are two partitions of $[n]$, then we define by $K\wedge
L$ the partition $\{K_i\cap L_j\}_{i=1}^k \vphantom{\}}_{j=1} ^l
\setminus \{\emptyset\}$.

\blemma\label{l3}
If the arguments of a function $f:{\Sigma}^n\to\Gamma$
are $K$- and $L$-separable, then they are $(K\wedge L)$-separable.
\elemma

\proofr Let
\begin{equation}\label{l3e}
f(\bar x)\equiv\bigoplus_{i=1}^k f_i(\bar x_{K_i})
 \equiv\bigoplus_{j=1}^l g_j(\bar x_{L_j}).
\end{equation}
For each $j$ from $1$ to $l$ we define the function $g'_j(\bar
x_{K_1\cap L_j},\ldots,\bar x_{K_k\cap L_j}) \triangleq g_i(\bar
x_{L_i})$, which differs from $g_j$ by the appropriate permutation of
the arguments, and the functions $h_{i,j}(\bar x_{K_i\cap
L_j})\triangleq g'_j(\bar 0,\ldots,\bar 0,\bar x_{K_i\cap L_j},\bar
0,\ldots,\bar 0)$. From (\ref{l3e}) for each $i$ from $1$ to $k$
we have
\begin{eqnarray*}
f_i(\bar x_{K_i}) \oplus \bigoplus_{i'\neq i} f_{i'}(\bar 0)   &
\equiv & \bigoplus_{j=1}^l g'_j(\bar 0,\ldots,\bar 0,\bar
x_{K_i\cap L_j},\bar 0,\ldots,\bar
0),\\
f_i(\bar x_{K_i}) & \equiv & \bigoplus_{j=1}^l h_{i,j}(\bar
x_{K_i\cap L_j}) \oplus const,\\
\mbox{and}\qquad
f(\bar x)&\equiv&\bigoplus_{i=1}^k \bigoplus_{j=1}^l h_{i,j}(\bar
x_{K_i\cap L_j})\oplus const.
\end{eqnarray*}

Since $K_i\cap L_j=\emptyset \Longrightarrow h_{i,j}(\bar
x_{K_i\cap L_j})=const$, the lemma is proved. \proofend

\blemma\label{l4} The decomposition
$$ f(\bar x)\equiv\bigoplus_{i=1}^k f_i(\bar x_{K_i}) $$
of a function $f:{\Sigma}^n\to\Gamma$ into functions $f_i$ with
non-separable arguments is unique up to constant summands in $f_i$.
\elemma

\proofr Let $\overline K$ be the set of partitions $K$ of
$[n]$ for which the function $f$ has $K$-separable
arguments. By Lemma~\ref{l3}, $(\overline K,\wedge)$ is a
semilattice. As we can see from the proof of Lemma~\ref{l3}, only
the least element of this semilattice corresponds to decomposition
into functions with non-separable arguments. By Lemma~\ref{l2}
these functions are unique up to constant summand.
\proofend



\ifx\href\undefined \newcommand\href[2]{#2} \fi\ifx\url\undefined
  \newcommand\url[1]{\href{#1}{#1}} \fi\ifx\bbljan\undefined
  \newcommand\bbljan{Jan} \fi\ifx\bblfeb\undefined \newcommand\bblfeb{Feb}
  \fi\ifx\bblmar\undefined \newcommand\bblmar{March} \fi\ifx\bblapr\undefined
  \newcommand\bblapr{Apr} \fi\ifx\bblmay\undefined \newcommand\bblmay{May}
  \fi\providecommand\bbljun{June}\ifx\bbljul\undefined \newcommand\bbljul{July}
  \fi\ifx\bblaug\undefined \newcommand\bblaug{Aug} \fi\ifx\bblsep\undefined
  \newcommand\bblsep{Sep} \fi\ifx\bbloct\undefined \newcommand\bbloct{Oct}
  \fi\ifx\bblnov\undefined \newcommand\bblnov{Nov} \fi\ifx\bbldec\undefined
  \newcommand\bbldec{Dec} \fi

\end{document}